\documentclass{article}%
\usepackage{amsmath}
\usepackage{amsfonts}
\usepackage{amssymb}
\usepackage{graphicx}%
\newtheorem{theorem}{Theorem}

\newtheorem{exercise}[theorem]{Exercise}

\begin{document}

\title{A Comment on: `On Some Contradictory Computations in Multi-Dimensional Mathematics'}
\author{{\footnotesize E. Capelas de Oliveira and} {\footnotesize W. A. Rodrigues,
Jr.}\\{\footnotesize Institute of Mathematics, Statistics and Scientific Computation
}\\{\footnotesize IMECC, UNICAMP, CP 6065}\\{\footnotesize 13083-859 Campinas SP, Brazil.}\\{\footnotesize e-mail: capelas@ime.unicamp.br walrod@ime.unicamp.br}}
\maketitle

\begin{abstract}
In this paper we analyze the status of some\ `unbelievable results' presented
in the paper `On Some Contradictory Computations in Multi-Dimensional
Mathematics' \cite{car} published in \textit{Nonlinear Analysis, }a journal
indexed in the Science Citation Index. Among some of the unbelievable results
`proved' in the paper we can find statements like that: (i) a rotation
$\mathcal{T}_{\theta}:\mathbb{R}^{2}\rightarrow\mathbb{R}^{2}$, $\theta\neq
n\pi/2$, is inconsistent with arithmetic, (ii) complex number theory is
inconsistent. Besides these `results' of mathematical nature \cite{car},
offers also a `proof' that Special Relativity is inconsistent. Now, we are
left with only two options (a) the results of \cite{car} are correct and in
this case we need a revolution in Mathematics (and also in Physics) or (b) the
paper is a potpourri of nonsense. We show that option (b) is the correct one.
All `proofs' appearing in \cite{car} are trivially wrong, being based on a
poor knowledge of advanced calculus notions. There are many examples (some of
them discussed in \cite{cawa,olro1,olro2,olro3,rod1}) of complete wrong papers
using nonsequitur Mathematics in the Physics literature. Taking into account
also that a paper like \cite{car} appeared in a Mathematics journal we think
that it is time for editors and referees of scientific journals to become more
careful in order to avoid the dissemination of nonsense.

\end{abstract}

\section{Introduction}

From time to time we give the following exercise to some of our students: find
mathematical errors and inconsistencies in articles appearing in scientific
journals, books or in the arXiv.\footnote{As result of this activity we
eventually write some notes which in some specific cases are sent for
publication. As example of this activity we quote the following papers:
\cite{cawa,olro1,olro2,olro3,rod1}.} Of course a paper with the title `On Some
Contradictory Computations in Multi-Dimensional Mathematics' called
immediately our attention as a potential one suggesting a possible exercise.
When we read that paper we first thought that the author was
joking\footnote{Unfortunately this was not the case, i.e., author was not
joking, as we discover reading some other papers signed by Carvalho which have
been quoted in \cite{car}.}, that he wrote it only to prove that many referees
indeed do not understand absolutely nothing about many of papers for which
they wrote reports\footnote{Something we also said, e.g, in
\cite{cawa,olro1,olro2,olro3,rod1}.}. Indeed, how a red herring was not
immediately activated when reviewing a paper that claims among other results
that :

\medskip

(i) \textquotedblleft multi-variable mathematics is inconsistent with
arithmetic ($1=0$) and also auto-contradictory as calculus is part of this
theory\textquotedblright,

(ii) \textquotedblleft A rotation $\mathcal{T}_{\theta}$, $\theta\neq n\pi/2$,
is inconsistent with arithmetic\textquotedblright,

(iii) \textquotedblleft Complex number theory is
inconsistent\textquotedblright,

(iv) \textquotedblleft Lorentz's transformation is contradictory unless $v=0$
(in which case the transformation is the identity) i.e., SRT is
trivial\textquotedblright.

\medskip

Below we show explicitly that all the above claims are based on a single
misconception, which result from the fact that author of \cite{car} forget
some crucial results of advanced calculus.

\section{Critical Analysis of \cite{car}}

The paper under review is divided in five sections. In the introduction it is
said that the permanent requirement of consistence is primordial for exact
science. So, if the claims (i-iv) above `proved' by author of \cite{car} were
true it is just the time to stop doing mathematics for a while and certainly
stop using it as presently known in any `exact' science. After some confuse
observations (including some ones concerning the theory of relativity) he
recalls the merge of algebra and geometry introduced by Descartes and says:

\textquotedblleft Thus contemporary mathematics adopts the notion of a change
of variables (see following section) as a tool to choose among all the
coordinate systems, the one which better simplifies the problem under study.
It application began to spread by the middle of the 18th century. What would
happen if two of these allowed choices did show up to be contradictory?
Clearly, Descartes' idea should be refined, and modern mathematics
reformulated. Observe that many new concepts (the negative and the complex
numbers) and theories (the infinitesimal calculus, linear algebra,
multi-variable calculus, non-Euclidean geometries, topology, etc.) were
conceived as a by-product of Descartes' idea. The notion of change of
variables stands as a fundamental concept to any of those
theories.\textquotedblright

After that he introduces in Section 2 his arguments for considering the change
of variables in $\mathbb{R}^{2}$ (different form the identity) to be inconsistent.

One of his thesis in that Section is formulate as follows. Consider the
standard statement:

\medskip

\textquotedblleft Any linear transformation $\mathcal{T}:\mathbb{R}%
^{2}\rightarrow\mathbb{R}^{2}$, $\mathcal{T}(x,y)=(\xi,\eta)$, given
by\footnote{The equations used in \cite{car} will be numbered here by its
number there with a prefix c.}%
\begin{equation}
\binom{\xi}{\eta}=\left[
\begin{array}
[c]{cc}%
a & b\\
c & d
\end{array}
\right]  \binom{x}{y} \tag{c2.1}%
\end{equation}
with $\det M=ad-bc\neq0$, is a valid transformation of variables (coordinates)
in $\mathbb{R}^{2}$\textquotedblright.

Well, Carvalho claims that such statement is \textquotedblleft obviously
meaningless in nature". He endorses his argument with an almost
incomprehensible example with barrels of wine, which shows clearly that he
does not know how to use mathematical theory in applications to problems of
the physical world\footnote{More on this below.}. Moreover, he continues his
discourse saying: \textquotedblleft Thus, the only relevant semantic content
of the term \textit{valid} used in this statement must be \textit{consistent
with mathematics, but not necessarily always consistent with nature}%
\textquotedblright\textit{. }Next he says that \textquotedblleft the assumed
validity of the property must be proved before definition is
cast\textquotedblright\ and, of course, all this naive discourse is only a
prelude for his `proof' in Section 3 of his paper of the inconsistencies in
the theory of changes of variables.

To `prove' his main thesis author first says that \textquotedblleft in general
a transformation from $\mathbb{R}^{2}$ to $\mathbb{R}$ of the form $u=ax+by$,
where $a$ and $b$ are constants that do not depend on parameters is consistent
with arithmetic as it represents the sum of two multiplicative tables. In the
change of variables Eq.(c2.1) it is required that two independent
transformations of this kind, namely
\begin{equation}
\left\{
\begin{array}
[c]{c}%
\xi=ax+by\\
\eta=cx+dy
\end{array}
\right.  \tag{c3.2}%
\end{equation}
be simultaneously satisfied\textquotedblright.

From this he says that from Eq.(c3.2) (or what is the same, his Eq.(c2.1)) it
follows that recalling that from Eq.(c3.2)
\begin{align}
x  &  =\eta/c-(d/c)y,\label{1}\\
y  &  =\eta/d-(c/d)x\nonumber
\end{align}
it follows that\footnote{Of course, the correct equations are: $\xi=\frac
{a}{c}\eta+(b-\frac{ad}{c})y$ and $\xi=\frac{a}{c}\eta+(a-\frac{bc}{d})x$.}:%

\begin{equation}%
\begin{tabular}
[c]{|c|}\hline
$\displaystyle \xi=\frac{a}{c}\eta+(b-\frac{d}{c})y,$\\\hline
\end{tabular}
\ \ \ \ \ \ \ \ \ \ \tag{c3.3}%
\end{equation}
and%

\begin{equation}%
\begin{tabular}
[c]{|c|}\hline
$\displaystyle \xi=\frac{b}{d}\eta+(a-\frac{c}{d})x.$\\\hline
\end{tabular}
\ \ \ \ \ \ \ \ \ \ \ \tag{c3.4}%
\end{equation}

Then, he says that \textit{synchronization} implies that the values of the
variables $\xi,\eta,x$ and $y$ do not change when we use them in Eq.(c2.1),
Eq.(c3.3) or Eq.(c3.4). Next he calculates $\xi_{x}=a$, $\xi_{y}=b$, $\eta
_{x}=c$ and $\eta_{y}=d$ and says that those equations \textquotedblleft hold
at any point $(x,y)$ in Eq.(c2.1)\textquotedblright. Next he says that
\textquotedblleft the implicit function theorem and the chain rule imply that
$\xi_{\eta}=\xi_{x}x_{\eta}$ and $\xi_{\eta}=\xi_{y}y_{\eta}$ also hold
everywhere in $\mathbb{R}^{2}$ (recall that $abcd\neq0$). Then he says that
from Eq.(c2.1) and Eq.(c3.3) it follows that $\xi_{\eta}=\xi_{x}x_{\eta}=a/c$
\ and $\xi_{\eta}=\xi_{y}y_{\eta}=b/d$ everywhere in $\mathbb{R}^{2}$. Finally
he concludes:

\medskip

\textquotedblleft Hence rationality requires that we must have, say,
$\xi_{\eta}(0,0)=a/c=b/d$ (the origin $(0,0)$ is the same in any system of
coordinates) so that arithmetic implies that $ad-bc=0$, contradicting the
standing hypothesis that $ad-bc\neq0$. Consequently, if we assume both the
validity of calculus and of the change of variables Eq.(c2.1), the mathematics
that follows from this assumption is contradictory when $abcd\neq
0$.\textquotedblright

\medskip

From this he establishes his theorem:

\textquotedblleft\textbf{Theorem 3.1}. If $ad-bc\neq0$ and $abcd\neq0$ then
the adoption of the validity of the concept of change of coordinates Eq.(2.1)
becomes inconsistent with the concept of partial derivative.\textquotedblright

From this point Carvalho deduces using arguments similar to the ones he
employed in his `proof' of his Theorem 3.1 that several other theories of
mathematics are inconsistents, e.g., he says that d'Alembert solutions of the
one-dimensional wave equation is inconsistent, that coordinate transformations
which represent rotations in $\mathbb{R}^{2}$ are inconsistent unless the
rotation angle is a of $n\pi/2$. With this last result, he establishes the corollary:

\textquotedblleft\textbf{Corollary 4.2}. Complex number theory is
inconsistent.\textquotedblright

And he did not stop here, he also `shows' in Section 5, that Lorentz
transformations are mathematically inconsistent, thus implying that Special
Relativity is inconsistent.

Well, my dear reader, at this point we think that if you have had a reasonable
course of advanced calculus you already realized the conceptual errors of
\cite{car}. Indeed, in any reasonable course of advanced calculus the
following exercise is usually given to students to verify if they grasped the
main concepts involved in the implicit function theorem (and the chain rule).

\begin{exercise}
\label{crux}Let
\begin{equation}
\xi=f_{1}(x,y),\qquad\eta=f_{2}(x,y) \label{1new}%
\end{equation}
with $f_{1},f_{2}:\mathbb{R}^{2}\supset U\rightarrow\mathbb{R}$ two $C^{1}$
continuous differentiable functions. Prove that the necessary and sufficient
condition for the existence of a functional relation between $\xi$ and $\eta$
of the form
\begin{equation}
\mathbf{F}(\xi,\eta)=0, \label{2}%
\end{equation}
is that
\[
\frac{\partial(\xi,\eta)}{\partial(x,y)}=0,
\]
where $\frac{\partial(\xi,\eta)}{\partial(x,y)}$ denotes the Jacobian of the
transformation given by \emph{Eq.(\ref{1new})}.
\end{exercise}

If you did not succeed solving the exercise, you may find a detailed solution
in Chapter 6 (Exercise 35) of Spiegel's book \cite{SPIEGEL}, and in that case
we suggest that you take the opportunity to review the basic concepts of calculus.

Here we use the result of the Exercise \ref{crux} to analyze Carvalho's
Theorem 3.1.

In this case, Eq.(c2.1) or Eq.(c3.1) defines the $C^{1}$ differentiable
functions $f_{1}$ and $f_{2}$. However, as the original meaning of Eq.(c2.1)
is the one of a transformation of variables in $\mathbb{R}^{2}$ we need to
impose that $\frac{\partial(\xi,\eta)}{\partial(x,y)}=ad-bc\neq0$. Then the
result of Exercise \ref{crux} says that it does not exist a functional
relation $\mathbf{F}(\xi,\eta)=0$, i.e., we cannot find a function $\mathbf{f}
:\mathbb{R\rightarrow R}$, $\xi\mapsto\eta=\mathbf{f}(\xi)$.

Thus under these conditions the calculations of \cite{car} leading to
$\xi_{\eta}=a/c$ and \ $\xi_{\eta}=b/d$ (reproduced above) are
\textit{meaningless, }for there is then nonsense in taking the derivative
$\xi_{\eta}$.

Note that if $ad-bc=0$ we can have according to the result of Exercise
\ref{crux} a functional relation $\mathbf{F}(\xi,\eta)=0$. Indeed if, e.g.,
(i) $c=\alpha a$, $d=\alpha b$, $\alpha\neq0$ we have $\xi=(a/c)\eta$, or (ii)
if $b=\beta a$, $d=\beta c$, $\beta\neq0$ we have again $\xi=(a/c)\eta$, as it
may be.

\section{ Relativity and \cite{car}}

The analysis presented above is, of course, enough to convince any reader that
\cite{car} is a very bad paper and should never be published in any scientific
journal, in particular a Mathematics journal that is indexed in the SCI.
However, we shall comment on another of Carvalho's statement. In Section 3 of
his paper he recalls that in the Theory of Relativity Lorentz transformations
play (as well known) a distinguished role. The nontrivial part of a special
Lorentz transformation, also called a boost in the $x$-direction can obviously
be written in the form of a linear transformation $\mathbb{R}^{2}%
\rightarrow\mathbb{R}^{2}$,%
\begin{equation}
\binom{\xi}{\eta}=\frac{1}{\sqrt{1-v^{2}/c^{2}}}\left[
\begin{array}
[c]{cc}%
1 & -v\\
-v/c^{2} & 1
\end{array}
\right]  \binom{x}{t}.\label{3}%
\end{equation}

He says that in Eq.(\ref{3}) $v$ denotes a `natural speed' with $0<v^{2}\leq
c^{2}$, $c$ denoting the speed of light. This transformation, as recalled in
\cite{car} has the same form of Eq.(2.1) with%
\begin{equation}
a=d=\frac{1}{\sqrt{1-v^{2}/c^{2}}},\text{ }b=\frac{-v}{\sqrt{1-v^{2}/c^{2}}%
},\text{ }c=-\frac{v}{c^{2}}\frac{1}{\sqrt{1-v^{2}/c^{2}}}. \label{4}%
\end{equation}

Now, in Section 3 of \cite{car} he calls the function $b$ of the variable $v$
by
\begin{equation}
a=:a(v)=\frac{-v}{\sqrt{1-v^{2}/c^{2}}}.\tag{c3.1}%
\end{equation}

Then he states: \textquotedblleft then $a$ is also a speed for each $v$ in its
range\textquotedblright. From that statement he thought that it is necessary
to impose that $a^{2}<c^{2}$ and then concludes that Eq.(c3.1) is
auto-contradictory when $c^{2}/2$ $<v^{2}<c^{2}$. He even says that
\textquotedblleft its numerator and denominator do not agree about their
mathematical and natural meanings\textquotedblright.

Of course, every reasonable Physics or Mathematics student that attended
lectures on the Theory of Relativity knows that the imposition that
$a^{2}<c^{2}$ is nonsense. Although $a$ has dimension of a velocity it is not
a `natural speed', and so we can have $a^{2}>c^{2}$. As well known, what the
theory of relativity forbids is that $v^{2}>c^{2}$.

So, all criticisms of \cite{car} concerning inconsistencies of  Theory of
Relativity are not valid, for our author besides having forgotten some
fundamental results of calculus, also does not know the simple
\textit{physical }meaning of the variables used in a theory of
Physics\footnote{We left as exercise to the reader to find additional errors
of Section 5 of \cite{car}.}. To endorse our statement we ask you to go to
Section 3 of \ \cite{car} where you can read: \textquotedblleft The most
common idea of change (or transformation) of variables occurs in the real line
$\mathbb{R}$ (arithmetic). It is of the form $y=ax$, where $a\neq0$ is a
constant, $x$ is the \textquotedblleft old\textquotedblright\ and $y$ is the
\textquotedblleft new\textquotedblright\ variable. Its effect is simply a
change of scale: the new unit is $a$ times the old one.\textquotedblright

Well, suppose that $x$ and $y$ denotes the linear measure of a rod in two
different units, \textit{meter} and \textit{centimeter}. It is clear that the
number $y$ is greater than the number $x$. This happens because the unit
called centimeter is $1/a=1/100$ times the unit called meter and not that the
unit called centimeter is $a=100$ greater than the unit called meter.

\section{Conclusions}

In this paper we analyzed the claims of \cite{car} that several theories of
Mathematics and also the Theory of Relativity are inconsistents. We showed
that all inferences of \cite{car} are the result of a simple fact: author of
\cite{car} does not know the implicit function theorem and, of course, is not
able to solve Exercise \ref{crux}. Also, it must be said that the publication
of completely wrong papers containing a potpourri of nonsense Mathematics (and
also Physics) is becoming more and more routine in `scientific' journals,
books and, of course, also in the arXiv. In our papers
\cite{cawa,olro1,olro2,olro3,rod1} we discussed several examples of nonsense
Mathematics. Our examples,\footnote{We have, of course, many and many others
examples, besides the ones we quoted.} we believe are enough to claim that it
is arrived the time for editors to choose better referees for their journals,
which at least must know advanced calculus.

\end{document}